\documentclass[conference]{IEEEtran}
\IEEEoverridecommandlockouts                              

\usepackage{epsf}
\usepackage{epsfig}
\usepackage{times}
\usepackage{amsmath,amssymb,amsfonts}
\usepackage{latexsym}
\usepackage{url}
\usepackage{mathrsfs}
\usepackage{algorithm}
\usepackage{caption}
\usepackage{datetime}
\usepackage{bm}
\usepackage{algorithm,algorithmic}
\usepackage{subfigure}
\usepackage{framed} 
\usepackage[framed]{ntheorem}
\usepackage{cite}
\usepackage{amstext}  
\usepackage{array}      
\usepackage{mathabx}
\usepackage{color}
\usepackage{textcomp}

\newframedtheorem{frm-definition}{Definition}

\newtheorem{theorem}{Theorem}

\newtheorem{lemma}{Lemma}

\newtheorem{remark}{Remark}
\newtheorem{assumption}{Assumption}
\newtheorem{proposition}{Proposition}

\newcommand{\hL}{\mathcal{L}}

\newcommand{\cD}{{\cal D}}

\newcommand{\cN}{{\cal N}}

\newcommand{\cZ}{{\cal Z}}

\def\BibTeX{{\rm B\kern-.05em{\sc i\kern-.025em b}\kern-.08em
    T\kern-.1667em\lower.7ex\hbox{E}\kern-.125emX}}
    
\begin{document}

\title{Stochastic Dual Algorithm for Voltage Regulation in Distribution Networks with Discrete Loads
\thanks{This work was supported by the U.S. Department of Energy under Contract No. DE-AC36-08GO28308 with the National Renewable Energy Laboratory; funded provided by the Grid Modernization Laboratory Consortium.}}

\author{Xinyang Zhou$^{*}$,
        Zhiyuan Liu$^{*}$,
        Emiliano Dall'Anese$^{\dagger}$,
        Lijun Chen$^{*}$\\
$^{*}$College of Engineering and Applied Science, University of Colorado, Boulder, USA\\
$^{\dagger}$National Renewable Energy Laboratory, Golden, USA\\
Emails: \emph{\{xinyang.zhou, zhiyuan.liu, lijun.chen\}@colorado.edu}, \emph{emiliano.dallanese@nrel.gov}}

\maketitle



\begin{abstract}
This paper considers power distribution networks with distributed energy resources and designs an incentive-based algorithm that allows the network operator and customers to pursue given operational and economic objectives while concurrently ensuring that voltages are within prescribed limits. Heterogeneous DERs with continuous and discrete control commands are considered. We address four major challenges:~discrete decision variables, non-convexity due to a Stackelberg game structure, unavailability of private information from customers, and asynchronous operation. 
Starting from a non-convex setting, we develop a distributed stochastic dual algorithm that solves a relaxed problem,
and prove that the proposed algorithm achieves the global optimal solution of the original problem on average. Feasible values for discrete decision variables are also recovered. Stability of the algorithm is analytically established and numerically corroborated.  
\end{abstract}


\section{Introduction}


Market-based mechanisms to control distributed energy assets have been recently developed with the objective of incentivizing customers to provide ancillary services to the grid while maximizing their own economic benefits; see, e.g.,~\cite{Subbarao13, Maharjan13,Tushar14} 
and pertinent references therein. However,  demand-response  and market-based problem formulations (e.g.,~\cite{Subbarao13,li2016market}) do not generally consider power flows in the distribution network (hence, they are oblivious to the voltage fluctuations that emerge from adjustments in the DER output power generated/consumed) and, oftentimes, stability is not analytically proven. The frameworks proposed in, e.g.,~\cite{li2011optimal,LinaCDC15} offer a way to account for the nonlinear power flows, but either their applicability is limited to a restricted class of network topologies or they consider only controllability of real powers.   

In our previous work~\cite{zhou2017pricing}, we  proposed an incentive-based algorithm that allows the network operator and customers to pursue given social welfare optimization while ensuring that voltage magnitudes are within the prescribed limits. We formulated a social-welfare maximization problem that captures a variety of optimization objectives, hardware constraints, and the nonlinear power-flow equations governing the physics of distribution systems; for the latter, a linear approximation of the nonlinear power-flow equations (see e.g.,~\cite{bolognani15,multiphaseArxiv}) 
is used to enable the development of a computationally-tractable optimal coordination method. 
To solve a well-defined but non-convex social-welfare maximization problem, we addressed two major challenges in \cite{zhou2017pricing}: 1)~We reformulated the non-convex social-welfare maximization problem into a convex problem, which is proved to be an exact convex relaxation; and
2)~We proposed an iterative distributed algorithm wherein customers and the network operator achieve consensus on a set of net real and reactive power that optimizes the objectives of both parties while ensuring that voltages are within the prescribed limits. In our design, customers are \emph{not} required to share private information regarding their cost functions and the operating region of their loads/generators with the  network operator. 

However, \cite{zhou2017pricing} only considers continuous decision variables (i.e., continuous DER commands), whereas in practice many appliances operate with discrete decision variables---e.g., capacity banks, thermostatically controlled loads (TCLs), and electric vehicles (EVs). In this work, we extend the incentive-based framework developed in \cite{zhou2017pricing} to include discrete decision variables. In literature of power system, discrete variables are dealt with either in a deterministic way \cite{bernstein2015design, liu1992discrete} 
 which usually generates suboptimal solutions, or in a stochastic way, e.g., \cite{kim2013scalable,macfie2010proposed} which often lacks rigorous analytical performance characterization.






In this paper we address discrete decision variables with a stochastic algorithm, and we establish analytical results for its convergence. Specifically, we first relax the discrete feasible sets to their convex hull. We then propose a distributed stochastic dual algorithm to solve the relaxed problem while recovering feasible power set points for discrete devices at every iteration, where two timescales  are considered for devices of different updating frequencies. Eventually, we prove that the proposed algorithm converges to a random variable whose mean value coincides with the optimal solution of the relaxed problem. 
We also characterize the variance of the resultant voltage due to the stochastic process of discrete power rate recovering, and design a robust implementation for  the voltage bounds accordingly. For completeness, notice that alternative ways to deal with discrete variables can be found in, e.g.,\cite{bernstein2015design} and \cite{kim2013scalable} (with other types of optimality claims, if any).

The rest of this paper is organized as follows. Section~\ref{sec:model} models the system and formulates the problem. Section~\ref{sec:reform} relaxes the problem, and proposes a distributed algorithm together with convergence analysis and performance characterization. Section~\ref{sec:simulation} provides numerical examples, and Section~\ref{sec:conclusion} concludes the paper.

\section{System Model and Problem Formulation}\label{sec:model}


\subsection{Network Model}
Consider a distribution network with $N+1$ nodes collected in the set $\cN \cup \{0\}$ with $\cN:=\{1, ..., N\}$. Let $p_i \in \mathbb{R}$ and $q_i \in \mathbb{R}$ denote the aggregated real and reactive power injections at node $i \in \cN$,  $V_i \in \mathbb{C}$  the phasor for the line-to-ground voltage  at node $i$, and define $v_i := |V_i|$. 
We consider a general linearized relationship between approximated voltage magnitudes $\hat{v}$ and injected powers as:
\vspace{-2mm}
\begin{align} 
	v\approx\hat{v}  = R p + X q + a , \label{eq:approximate} 
\end{align}
where matrices $R, X \in \mathbb{R}^{N \times N}$ and the vector $a \in \mathbb{R}^{N}$ can be formed as described in, e.g.,~\cite{bolognani15,multiphaseArxiv,farivar2013equilibrium}.

\begin{remark}
The linear model \eqref{eq:approximate} is utilized to facilitate the design of computationally-affordable algorithms. Such approximation usually introduces a bounded error in voltages (see e.g., \cite{chiang1990existence, zhou2016vvac}), 
which further results in a bounded discrepancy between the generated and the optimal operating points. In Section~\ref{sec:simulation}, we will however utilize nonlinear power flow model for numerical examples. 
\end{remark}

\subsection{Node and Device Model}
\subsubsection{Nodal power aggregation}
At node $i\in\cN$, assume that the aggregated power is from two kinds of devices: non-controllable  and controllable. Denote by $p_{i,0}\in \mathbb{R}$ and $q_{i,0}\in \mathbb{R}$ the total real and reactive power from non-controllable devices. Meanwhile, assume a customer $i$ controls the devices collected in a set $\cD_i$. Denote by $p_{i,d} \in \mathbb{R}$ and $q_{i,d} \in \mathbb{R}$ the real and reactive power injections of controllable device $d\in\cD_i$. Then the power aggregation at node $i$ is cast as $p_{i}=p_{i,0}+\sum_{d\in\cD_i}p_{i,d}$ and $q_{i} =q_{i,0}+\sum_{d\in\cD_i}q_{i,d}$.

\subsubsection{Devices}
At node $i$, we consider the following two types of controllable devices: 
(i) continuous-rate fast-updating devices (\textbf{fast devices}), e.g., PV inverters, collected in a subset $\cD_{F_i}\subseteq\cD_i$, and
(ii) discrete-rate slow-updating devices (\textbf{slow devices}), e.g., TCLs, collected in a subset $\cD_{S_i}\subseteq\cD_i$.
We also assume that the slow devices only consume real power without loosing generality.

	Let $\mathcal{Z}_{F_i,d}$ denote  the  feasible power injection set of fast devices $d\in\cD_{F_i}$, e.g., $\mathcal{Z}_{F_i,d}$ of a PV system has the form of:
\begin{align} 
\cZ_{F_i,d} =  \left\{(p_{i,d}, q_{i,d}) \hspace{-.1cm} :  0 \leq p_{i,d}  \leq  p_{i,d}^{\text{av}}, p_{i,d}^2 + q_{i,d}^2 \leq  \eta_{i,d}^2 \right\},\nonumber
\end{align}
where $p_{i,d}^{\text{av}}$ denotes the available real power from a PV system, and $\eta_{i,d}$ is the rated apparent power capacity. 
	
	Meanwhile, considering that the slow devices $d\in\cD_{S_i}$ have discrete power rates, and possibly some temporal dynamics and constraints that can be modeled as affine equalities and inequalities (e.g., temperature requirement of TCLs, state of charge of battery, etc.), we denote by $\mathcal{Z}_{S_i,d}$ the \emph{convex hull} of the feasible sets of slow devices. For example, an air conditioner follows the temperature dynamics with constraints:
	\begin{eqnarray}
	&&T_{i,d}^{\text{in}+}=T^{\text{in}}_{i,d}+\theta_1(T^{\text{out}}_{i,d}-T^{\text{in}}_{i,d})+\theta_2 p_{i,d},\label{eq:tcl1}\\
	&&\underline{T}_{i,d}^{\text{in}}\leq T_{i,d}^{\text{in}+}\leq\overline{T}_{i,d}^{\text{in}},\label{eq:tcl2}
	\end{eqnarray}
	where the room temperature in the next considered moment (e.g., 15 minutes ahead) $T_{i,d}^{\text{in}+}$ is a linear function of power consumption rate $p_{i,d}$ given current room temperature $T^{\text{in}}_{i,d}$ and ambient temperature $T^{\text{out}}_{i,d}$ with constants $\theta_1$ and $\theta_2$ related to thermal properties and time interval length, and $\underline{T}_{i,d}^{\text{in}}$ and $\overline{T}_{i,d}^{\text{in}}$ being customized temperature bounds.
Now assume the air conditioner has three power consumption rates $0< p_{\text{low}}< p_{\text{high}}$. Then its set $\mathcal{Z}_{S_i,d}$ is computed as: 
\begin{align} 
\cZ_{S_i,d} =\{p_{i,d}|(\ref{eq:tcl1})\text{--}(\ref{eq:tcl2})~\text{and}~p_{i,d}\in[0,p_{\text{high}}]\}.\nonumber
\end{align}
	
	 Overall, $\mathcal{Z}_{F_i,d}$ and $\mathcal{Z}_{S_i,d}$ are assumed to be convex and compact, and so are the following sets (the symbol $\bigtimes$ denotes the Cartesian product):
	 \begin{align*}
	 &\cZ_{F_i}:=\underset{d\in\cD_{F_i}}{\bigtimes}\cZ_{F_i,d},&&\hspace{-8mm}\cZ_{S_i}:=\underset{d\in\cD_{S_i}}{\bigtimes}\cZ_{S_i,d};\\[-3pt]
	 &\cZ_F:=\underset{i\in\cN}{\bigtimes}\cZ_{F_i},&&\hspace{-8mm}\cZ_S:=\underset{i\in\cN}{\bigtimes}\cZ_{S_i};\\[-3pt]
 	 &\cZ_i\!\!:=\!\!\mathcal{Z}_{F_i}\bigtimes\mathcal{Z}_{S_i},&&\hspace{-8mm}\cZ:=\cZ_{F} {\bigtimes}\cZ_{S}.
	 \end{align*}
	Further, for simplicity, we introduce the following notation:
	\begin{align*}
	&z_{_{F_i}}= \{p_{i,d},q_{i,d}\}_{d\in\cD_{F_i}}\!\in\!\cZ_{F_i},&& z_{_{S_i}}= \{p_{i,d},q_{i,d}\}_{d\in\cD_{S_i}}\!\in\!\cZ_{S_i};\\
	&z_{_F}= \{z_{_{F_i}}\}_{i\in\cN}\!\in\!\cZ_F,&& z_{_S}= \{z_{_{S_i}}\}_{i\in\cN}\!\in\!\cZ_S;\\
	&z_i= \{p_{i,d},q_{i,d}\}_{d\in\cD_i}\!\in\!\cZ_i, && z = \{z_i\}_{i\in\cN}\!\in\!\cZ.
	\end{align*}
	Consider a cost function $C_{i}(z_i)$ that captures well-defined performance objectives of all devices at node $i$. The following assumption is made.
	\begin{assumption}
		\label{as:Cinv}
		Functions $C_{i}(z_i),\ \forall i\in\cN$ are continuously differentiable and strongly convex in $z_i$. Moreover, the first-order derivative of $C_{i}(z_i)$ is bounded in $\mathcal{Z}_{i}$.
	\end{assumption}
%

\subsection{Problem Formulation}

The overall objective is to design a strategy where the net- work operator and the customers pursue their own operational goals and economic objectives, while also achieving global coordination to enforce voltage regulation. In our design, incentive signals are designed to achieve this goal.

\vspace{.1cm}

\subsubsection{Customer's problem} 
Let  $\alpha_i\in\mathbb{R}$ and $\beta_i\in\mathbb{R}$ be the pricing/reward signals sent by the network operator (e.g., distribution system operator or aggregator) to customer $i$ for  real and reactive power injections, respectively. Customers are assumed to be rational and self-interested, aiming to minimize their own cost, by solving the following problem $(\bm{\mathcal{P}_{1,i}})$ given signals $(\alpha_i,\beta_i)$:
\vspace{-.1cm}
\begin{subequations}
	\begin{eqnarray}
	&\hspace{-6mm}\underset{z_i}{\min}&\!\!\!C_{i}(z_i)-\sum_{d\in\cD_i}(\alpha_i p_{i,d}+\beta_i q_{i,d}),\label{eq:obj_c}\\[-4pt]
	&\hspace{-6mm}\mathrm{s.t.}&z_i \in \mathcal{Z}_{i}, \label{eq:pccon_c}
	\end{eqnarray}
\end{subequations}
where $\alpha_i p_{i,d}$ and $\beta_i q_{i,d}$ represent payment to/reward from the network operator owing to device $d\in\cD_i$.

Because (\ref{eq:obj_c}) is strongly convex in $z_i$, and $\mathcal{Z}_{i}$ is convex and compact,  a unique solution  $z_i^*$ exists. 
For future developments, consider a ``best-response" strategy $b_{i}$ of customer $i$ as the following function of $(\alpha_i,\beta_i)$:
\begin{eqnarray}
z_i^*=b_{i}(\alpha_i,\beta_i):=\underset{z_i \in \mathcal{Z}_{i}}{\arg\min}~ C_{i}(z_i)-\sum_{d\in\cD_i}(\alpha_i p_{i,d}+\beta_i q_{i,d}).\label{eq:br}\nonumber
\end{eqnarray}

\subsubsection{Recover feasible power rates}
Given $p^{*}_{i,d}$ solved from the relaxed feasible set for slow device $d$, we randomly select a feasible power rate $p_{i,d}$ based on the probability distribution such that  $E[p_{i,d}]=p_{i,d}^{*}$,
where $E[ \cdot]$ denotes the expectation.

While there are multiple ways to determine the probability distribution of feasible setpoints based on  $p^{*}_{i,d}$, we exemplify the procedure with two-point distribution for illustrative purpose. We select two feasible power rates and denote them as $\underline{p}_{i,d}$ and $\overline{p}_{i,d}$, such that $\underline{p}_{i,d}\leq p^{*}_{i,d}\leq \overline{p}_{i,d}$. 
Then the related probability of the corresponding two-point distribution is calculated as:
\begin{eqnarray}\label{eq:probability}
	\left \{ 
  	\begin{array}{l l l }
	\mathbb{P}(p^t_{i,d}=\underline{p}_{i,d})&=&({\overline{p}_{i,d}-p_{i,d}^{t*}})/({\overline{p}_{i,d}-\underline{p}_{i,d}}),\\
	\mathbb{P}(p^t_{i,d}=\overline{p}_{i,d})&=&({p_{i,d}^{t*}-\underline{p}_{i,d}})/({\overline{p}_{i,d}-\underline{p}_{i,d}}),
	\end{array}\right.
\end{eqnarray}
according to which $p_{i,d}$ is randomly chosen.

\subsubsection{Social-welfare problem} 
Consider the following optimization problem $(\bm{\mathcal{P}_2})$ capturing both social welfare and voltage constraints of a distribution feeder:
\begin{subequations}
\begin{eqnarray}
 &\hspace{-6mm}\underset{z,\hat{v}, \alpha,\beta}{\min}& \sum_{i\in\cN}C_{i}(z_i),\label{eq:obj}\\[-4pt]
&\hspace{-6mm}\mathrm{s.t.} & p_i=p_i^0+\!\!\sum_{d\in\cD_i}p_{i,d},\ q_i=q_i^0+\!\!\sum_{d\in\cD_i}q_{i,d},\\[-2pt]
&\hspace{-6mm}&\hat{v}= Rp + Xq + a, \label{eq:volt0}\\
&\hspace{-6mm}&\underline{v}\leq \hat{v} \leq \overline{v},\label{eq:volt}\\
&\hspace{-6mm}& z_i=b_{i}(\alpha_i,\beta_i),i\in\cN,\label{eq:pccon}
\end{eqnarray}
\end{subequations}
where vectors $\underline{v}$ and $\overline{v}$ are prescribed minimum and maximum voltage magnitude limits (e.g., ANSI C84.1 limits) enforced by the network operator. Notice that the total payment from/reward to the customers cancels out the total payment received/reward paid by the network operator; thus, it is not in the social welfare objective function (\ref{eq:obj}). Notice that the problem formulation naturally extends to the case where multiple customers are located at node $i$,  but we outline the problem in this way to limit the complexity of the notation.

The best-response strategy of $(\bm{\mathcal{P}_{1,i}})$ is embedded in $(\bm{\mathcal{P}_2})$ through the constraint (\ref{eq:pccon})---i.e., the network operator has knowledge of the reaction of the customers toward any signals and takes it into consideration when making decisions. This constitutes a Stackelberg game wherein the network operator solves $(\bm{\mathcal{P}_2})$ and sends the signals $(\alpha^*,\beta^*)$ to customers;  subsequently, each customer responds with computed power injections $z_i^*$ from $(\bm{\mathcal{P}_{1,i}})$ based on the received signals. By design, $z_i^*$ coincides with the optimal solution of $(\bm{\mathcal{P}_2})$. 

However, it is impractical to solve problem $(\bm{\mathcal{P}_2})$ not only because of its non-convexity introduced by constraints (\ref{eq:pccon}), but also because of its requirement for customers' full information, which is usually private.  In the following section, we first reformulate $(\bm{\mathcal{P}_2})$ as a convex optimization problem with a signal design strategy; they together bypass the problem of non-convexity, while achieving the global optimal solution of $(\bm{\mathcal{P}_2})$. Then, based on the stochastic dual algorithm, we propose a distributed algorithm that prevents any exposure of private information from the customers while solving a  convex optimization problem where  devices admits  discrete power levels.

\section{Distributed Algorithm Design}\label{sec:reform}
          
\subsection{Convex Reformulation and Signal Design}
Consider the following  convex optimization problem $(\bm{\mathcal{P}_3})$:
\vspace{-5mm}
\begin{subequations}
\begin{eqnarray}
 &\hspace{-6mm}\underset{z, \hat{v}}{\min}& \sum_{i\in\cN}C_{i}(z_i),\label{eq:obj2}\\[-4pt]
&\hspace{-6mm}\mathrm{s.t.} & p_i=p_i^0+\!\!\sum_{d\in\cD_i}\!p_{i,d},\ q_i=q_i^0+\!\!\sum_{d\in\cD_i}\!q_{i,d},\\[-2pt]
&\hspace{-6mm}&\hat{v}= Rp + Xq + a, \label{eq:volt02}\\
&\hspace{-6mm}&\underline{v}\leq \hat{v} \leq \overline{v},\label{eq:volt2}\\
&\hspace{-6mm}& z_i\in \mathcal{Z}_{i}, \forall i\in\cN, \label{eq:pccon2}
\end{eqnarray}
\end{subequations}
with non-convex constraint (\ref{eq:pccon}) replaced with (\ref{eq:pccon2}), and $\alpha$, $\beta$ to be determined later. We assume $(\bm{\mathcal{P}_3})$ is feasible:

\begin{assumption}
\label{ass:sla}
Slater's condition holds for $(\bm{\mathcal{P}_3})$. 

\end{assumption}

Given the strong convexity of the objective function \eqref{eq:obj2} in $z$, a unique optimal solution exists for problem $(\bm{\mathcal{P}_3})$. Notice that a solution $(z^*, {\hat{v}}^{*})$ of $(\bm{\mathcal{P}_3})$ might not be feasible for $(\bm{\mathcal{P}_2})$, because there might not exist a $(\alpha^*, \beta^*)$ such that $z^*_i=b_i (\alpha_i^*,\beta_i^*)$. We will, however, show next that such $(\alpha^*, \beta^*)$ exists; thus, the solution of $(\bm{\mathcal{P}_3})$ gives the solution of $(\bm{\mathcal{P}_2})$.

Substitute (\ref{eq:volt02}) into (\ref{eq:volt2}), and denote by $\underline{\mu}$ and $\overline{\mu}$ the dual variables  associated with the constraints~\eqref{eq:volt2}. Let $\hat{v}^*$ be the optimal voltage magnitudes produced by $(\bm{\mathcal{P}_3})$, and denote the optimal dual variables associated with~\eqref{eq:volt2} as $\underline{\mu}^*,\overline{\mu}^*$. Then, we design the signals as:
	\begin{eqnarray}
		\alpha^*= R\big[\underline{\mu}^*-\overline{\mu}^*\big],\ \ \beta^*= X\big[\underline{\mu}^*-\overline{\mu}^*\big].\label{eq:signal}
	\end{eqnarray}
Note that $\alpha^*$ and $\beta^*$ are composed of dual prices $\underline{\mu}^*$ and $ \overline{\mu}^*$ with $R, X$ characterizing the network structure. We can prove that the above signals are bounded, precluding the possibility of infinitely large signals. 

\begin{theorem}[Theorem 1 of \cite{zhou2017pricing}]\label{lem:boundmu}
Under Assumptions~\ref{as:Cinv}--\ref{ass:sla}, the signals $(\alpha^*, \beta^*)$ defined by \eqref{eq:signal}  are bounded.  
\end{theorem}


By examining the optimality conditions of $(\bm{\mathcal{P}_2})$ and $(\bm{\mathcal{P}_3})$, we have the following result.

\begin{theorem}[Theorem 2 of \cite{zhou2017pricing}]
\label{the1}
The solution of problem $(\bm{\mathcal{P}_3})$ along with the signals $(\alpha^*, \beta^*)$ defined in (\ref{eq:signal}) is a global optimal solution of problem $(\bm{\mathcal{P}_2})$; i.e., problem  $(\bm{\mathcal{P}_3})$  is an exact convex relaxation of problem  $(\bm{\mathcal{P}_2})$. 
\end{theorem}
 
\vspace{.1cm}

From now on, we will use the optima of $(\bm{\mathcal{P}_3})$ and $(\bm{\mathcal{P}_2})$ interchangeably depending on the context. Next, based on Theorem \ref{the1}, we will develop an iterative algorithm that achieves the optimum of $(\bm{\mathcal{P}_3})$ (and hence that of $(\bm{\mathcal{P}_2})$) without exposing any private information of the customers to the network operator.


\subsection{Two Timescales and Iterative Algorithm}\label{sec:iter}
In this part, we design an iterative algorithm to solve $(\bm{\mathcal{P}_3})$. As mentioned, we have two types of devices with two different update frequencies. Assume that slow devices update $M$ times slower than fast devices with integer $M\geq 1$. We index by $k\in\mathbb{Z}_{++}$ the iterations when fast devices update. Then slow devices updates when $k=tM$ with index $t\in\mathbb{Z}_{++}$. We put the two timescales update in Algorithm~\ref{alg:twotime} for easy reference later. Based on this  strategy, we next propose a stochastic dual algorithm to solve $(\bm{\mathcal{P}_3})$ while recovering feasible power rates for slow devices.

	\begin{algorithm}
		\caption{Two-timescale update} \label{alg:twotime}
		\begin{algorithmic}
			\IF{iteration $k=tM$}
			\STATE Customer $i$ solves $z_i^*(k+1)$  
			\begin{eqnarray}
			=\underset{ z_i \in \mathcal{Z}_{i}}{\arg\min}~ C_{i}(z_i)-\!\!\sum_{d\in\cD_i}\big(\alpha_i(k) p_{i,d}+\beta_i(k) q_{i,d}\big),\nonumber		 	
			\end{eqnarray}
			recovers $z_{i,d}(k+1)$ with $E[z_{i,d}(k+1)]=z_{i,d}^*(k+1)$ for $d\in\cD_{S_i}$, and sets $z_{i,d}(k+1)=z_{i,d}^*(k+1)$ for $d\in\cD_{F_i}$.
			\ELSIF{iteration $k=tM+m,\ m=1,\ldots,M-1$}
			\STATE Customer $i$ keeps $z_{i,d}(k+1)=z_{i,d}(k)$ for $d\in\cD_{S_i}$, and gets $z_{i,d}(k+1)$ for $d\in\cD_{F_i}$ by solving:
			\begin{eqnarray}
			\underset{ z_{_{F_i}} \in \mathcal{Z}_{F_i}}{\arg\min}~ C_{i}(z_{_{F_i}} | z_{_{S_i}})-\!\!\sum_{d\in\cD_{F_i}}(\alpha_i(k) p_{i,d}+\beta_i(k) q_{i,d}),\nonumber
			\end{eqnarray}
			\vspace{-.3cm}
			\ENDIF
		\end{algorithmic}
	\end{algorithm}

Denote by $\mu:=[\underline{\mu}^{\intercal}, \overline{\mu}^{\intercal}]^{\intercal}\in\mathbb{R}^{2N}_+$ the vector of stacked dual variables, and denote $g(z)=\begin{bmatrix}\underline{v}-Rp-Xq-a\\Rp+Xq+a-\overline{v} \end{bmatrix}$. We can write the Lagrangian of $(\bm{\mathcal{P}_3})$ as:
\vspace{-.0cm}
\begin{eqnarray}
\hL(z,\mu)=\hL(z_{_S},z_{_F},\mu)=\underset{i\in\cN}{\sum} C_{i}(z_i)+\mu^{\intercal} g(z),  \label{eq:Lag}
\end{eqnarray}
which is obtained by keeping the constraints $z\in\mathcal{Z}$ and $\mu\in\mathbb{R}_+^{2N}$ implicit. 
Fix the value of $z_{_S}$, and define the resultant form as a ``reduced" Lagrangian $\hL_F(z_{_F},\mu|z_{_S})$.

We will implement a dual algorithm with two timescales to solve the following minimax problems based on $\hL(z,\mu)$ and $\hL_F(z_{_F},\mu | z_{_S})$:
\begin{eqnarray}
\underset{\mu\in\mathbb{R}_+^{2N}}{\max} \underset{z\in\mathcal{Z}}{\min} ~\hL(z,\mu),\ \text{and}\ \ \underset{\mu\in\mathbb{R}_+^{2N}}{\max} \underset{z_{_F}\in\mathcal{Z}_F}{\min} \hL(z_{_F},\mu|z_{_S}).\label{eq:minmax}
\end{eqnarray}

To this end, we define two concave dual functions for $\hL$ and $\hL_F$, respectively: 
\vspace{-.1cm}
\begin{eqnarray}
h(\mu):=\underset{z\in\mathcal{Z}}{\min}~\hL(z,\mu),\ \text{and}\ \ h_F(\mu|z_{_S}):=\underset{z_{_F}\in\mathcal{Z}_F}{\min}\hL_F(z_{_F},\mu|z_{_S}),\nonumber
\end{eqnarray}
with corresponding dual problems:
\vspace{-.1cm}
\begin{eqnarray}
\underset{\mu\in\mathbb{R}_+^{2N}}{\max}~h(\mu), \ \text{and}\ \ \underset{\mu\in\mathbb{R}_+^{2N}}{\max}~h_F(\mu|z_{_S}). \nonumber
\end{eqnarray}

Considering the dual algorithm to solve (\ref{eq:minmax}) while recovering implementable feasible power rates for discrete devices at each iteration, we have the following stochastic dual algorithm:
\vspace{-.5cm}
\begin{subequations}\label{eq:dualalg}
\begin{eqnarray}
z(k+1)&&\text{set by Algorithm~\ref{alg:twotime}},\label{eq:random}\\
\mu(k+1)&=&\big[\mu(k)+\varepsilon_k g(z(k+1))\big]_+,
\end{eqnarray}
\end{subequations}
where $[~]_+$ is a projection operator onto the nonnegative orthant, and the stepsize $\varepsilon_k$ (we will show how to select $\varepsilon_k$ shortly). Also notice that $g(z(k+1))$ is subgradient of $h(\mu(k))$ when $k=tM$ and that of $h_F(\mu(k))$ when $k=tM+m,\ m=1,\ldots,M-1$.

Because $\mathcal{Z}$ is compact and $g(z)$ is linear in $z$, there exists some constant $G>0$ such that $E[\|g(z(k))\|]\leq G$ for all $k$. We next show the stability of dynamics~(\ref{eq:dualalg}) with diminishing stepsize in the following subsection.

\subsection{Performance Analysis with Diminishing Stepsize}\label{sec:diminish}
In this part, we choose stepsize $\varepsilon_k$ to be square summable but not summable, i.e.:
\vspace{-.0cm}
\begin{eqnarray}
	\sum_{k=1}^{\infty}\varepsilon_k^2 <\infty,\ \ \sum_{k=1}^{\infty}\varepsilon_k =\infty,\label{eq:stepsize}
\end{eqnarray}
e.g., $\varepsilon_k=1/t$ at iteration $k=tM+m,\ m=0,\ldots, M-1$. With such diminishing stepsize, we will prove the convergence of the sequence $\{\mu_k\}$ generated by (\ref{eq:dualalg}) to a random vector. Moreover, we characterize the variance of voltage due to the randomness, and we propose a robust design.

\subsubsection{Convergence}
To show the convergence of dynamics~(\ref{eq:dualalg}), we utilize the next lemma \cite{polyakintroduction}.

\begin{lemma}\label{gladyshev}
Consider a sequence of random variables $\omega(1),\ldots, \omega(k)\geq 0$, $E[\omega(1)]\!<\!\infty$ and
$E[\omega(k+1)|\omega(1),\ldots, \omega(k)]\!\leq\! (1\!+\!x_k)\omega(k)\!+\!y_k$,
with $\sum_{k=1}^{\infty}x_k\!<\!\infty,\ \sum_{k=1}^{\infty}y_k\!<\!\infty,\ x_k\!\geq\! 0,\ y_k\!\geq\! 0$. Then $\omega(k)\rightarrow \omega(\infty)$ almost surely, where $\omega(\infty)\geq 0$ is some random variable.
\end{lemma}

We then have the following convergence result, the proof of which is referred to the Appendix.

\begin{theorem}\label{the:convergence}
If the stepsize $\varepsilon_k$ is chosen as in (\ref{eq:stepsize}), the sequence $\{\mu(k)\}$ generated by~(\ref{eq:dualalg}) converges to certain random vector $\mu(\infty)$ almost surely.
\end{theorem}

Denote $\tilde{\mu}(k):=E[\mu(k)]$.
By Theorem \ref{the:convergence}, $\lim_{k\rightarrow\infty}\tilde{\mu}(k)=\tilde{\mu}(\infty)$, where $\tilde{\mu}(\infty)$ is the mean value of random variable ${\mu}(\infty)$ to which $\{\mu(k)\}$ converges.
We next show that $\tilde{\mu}(\infty)=\mu^*$. The proof is also referred to the Appendix. 

\begin{theorem}\label{the:performance}
	Select the stepsize $\varepsilon_k$ as in (\ref{eq:stepsize}). The sequence $\{\tilde{\mu}(k)\}$ generated by~(\ref{eq:dualalg}) converges to $\mu^*$. Meanwhile, the running average of $h(\mu(k))$ approaches $h(\mu^*)$ as $k\rightarrow\infty$; i.e.:
	\vspace{-.1cm}
	\begin{eqnarray}
	&&\lim_{k\rightarrow\infty} h(\mu^*)-\sum_{\kappa=1}^k {h(\mu(\kappa))}/{k}=0.\label{eq:timeaverage} 
	\end{eqnarray}	
\end{theorem}

\begin{remark}
With strongly convex cost functions $C_i$, $\{\mu(k)\}$ generated by~(\ref{eq:dualalg}) usually converges to a random vector even with constant stepsize, where all the properties we obtain in Theorem~\ref{the:convergence}--\ref{the:performance}
hold. This will be shown with numerical examples in Section~\ref{sec:scenarios}.
\end{remark}

\subsubsection{Variance and Robust Design}
The randomness in the power rate selection for discrete devices of $\cD_{S}$ leads to volatility of voltages. Let $D_{S}$ be the number of all discrete devices. We next characterize the variance of the voltage (the result is tailored to the variable \eqref{eq:probability}, but can be straightforwardly generalized). 

\begin{proposition}\label{pro:var}
Using the randomized selection strategy \eqref{eq:probability}, the voltage variance $Var(v_i)$ at node $i\in\cN$ is bounded as:
\begin{eqnarray}
Var(\hat{v}_i)\leq D_{S}/4\sum_{j\in\cN}R^2_{ij}\cdot \underset{j,d}{\max} (\overline{p}_{j,d}-\underline{p}_{j,d})^2.\label{eq:variance}
\end{eqnarray}
\end{proposition}
The proof is referred to the Appendix.

The result of Proposition~\ref{pro:var}  motivates us to propose the following robust design. We choose tighter voltage bounds $[\underline{v}',\overline{v}']$ with $\underline{v}<\underline{v}'<v^{\text{nom}}<\overline{v}'<\overline{v}$, and replace the original bounds with the tighter bounds in the algorithm so that the resultant voltage falls within the original bounds with required probability. We design the tighter bounds based on Chebyshev's inequality \cite{stark2014probability} as shown next. 




\begin{proposition}[robust implementation]\label{pro:rob}
Given $\delta>0$, if the  voltage upper and lower bounds are set as:
$\overline{v}'_i \leq \overline{v}_i-\delta$ and $\underline{v}'_i \geq \underline{v}_i+\delta$,
then:
\begin{eqnarray}
\mathbb{P}[\hat{v}_i\geq \overline{v}_i]\leq{Var(\hat{v}_i)}/{2\delta^2},\ \text{and}\ \ \mathbb{P}[\hat{v}_i\leq \underline{v}_i]\leq{Var(\hat{v}_i)}/{2\delta^2}.\nonumber
\end{eqnarray}
\end{proposition}
Please see the proof in the Appendix.
\begin{remark}
	These bounds are admittedly conservative; however, reasonable values can be obtained in realistic settings. For example, scenario 2) in Section \ref{sec:discretesetup} with variance estimated by upper-bound \eqref{eq:variance} leads to robust bounds $\overline{v}'_i=1.035\!$ p.u. and $\underline{v}'_i=0.965\!$ p.u., with ${Var(\hat{v}_i)}/{2\delta^2}\leq 5\%$ and $\overline{v}_i=1.05\!$ p.u. and $\underline{v}_i=0.95\!$ p.u. Tighter bounds can be obtained empirically.
\end{remark}

%
%

\subsection{Distributed Stochastic Dual Algorithm}
The decomposable structure of (\ref{eq:dualalg}) naturally enables the following iterative \emph{distributed} algorithm:
\vspace{-.1cm}
\begin{subequations}\label{eq:dyn2}
		\begin{eqnarray}
		\hspace{-4mm}z(k+1)&& \text{set by Algorithm~\ref{alg:twotime}},\label{eq:dyncus2} \\
		\hspace{-4mm}\hat{v}(k+1)&=& Rp(k+1)+Xq(k+1) + a \label{eq:dynv}\\
		\hspace{-4mm}\underline{\mu}(k+1)&=& \big[\underline{\mu}(k)+\varepsilon_k\big(\underline{v}-\hat{v}(k+1)\big)\big]_+\,,\label{eq:dynmu}\\
		\hspace{-4mm}\overline{\mu}(k+1)&=& \big[\overline{\mu}(k)+\varepsilon_k\big(\hat{v}(k+1)-\overline{v}\big)\big]_+\,,\\
		\hspace{-4mm}\alpha(k+1)&=&R\big[\underline{\mu}(k+1)-\overline{\mu}(k+1)\big]\,,\label{eq:dynalpha}\\
		\hspace{-4mm}\beta(k+1)&=& X \big[\underline{\mu}(k+1)-\overline{\mu}(k+1)\big]\,,\label{eq:dynbeta}
		\end{eqnarray}
\end{subequations}
where the power set points of the devices are computed and implemented locally through~\eqref{eq:dyncus2}, and~\eqref{eq:dynv}--\eqref{eq:dynbeta} are performed at the network operator. Notice that each customer $i$ \emph{does not} share its cost function $C_{i}$ or its feasible set $\mathcal{Z}_{i}$ with the network operator; rather, the customer transmits to the network operator only the resultant power injections $z_{i}(k)$. Indeed,(\ref{eq:dualalg}) and (\ref{eq:dyn2}) are equivalent and the results of Theorem \ref{the:convergence}--\ref{the:performance} apply to (\ref{eq:dyn2}).

%
%
%
%
%
%
%


\section{Application Scenarios}\label{sec:simulation}
\label{sec:scenarios}
\begin{figure}
	\centering
	\vspace{.25cm}
	\includegraphics[width=.4\textwidth]{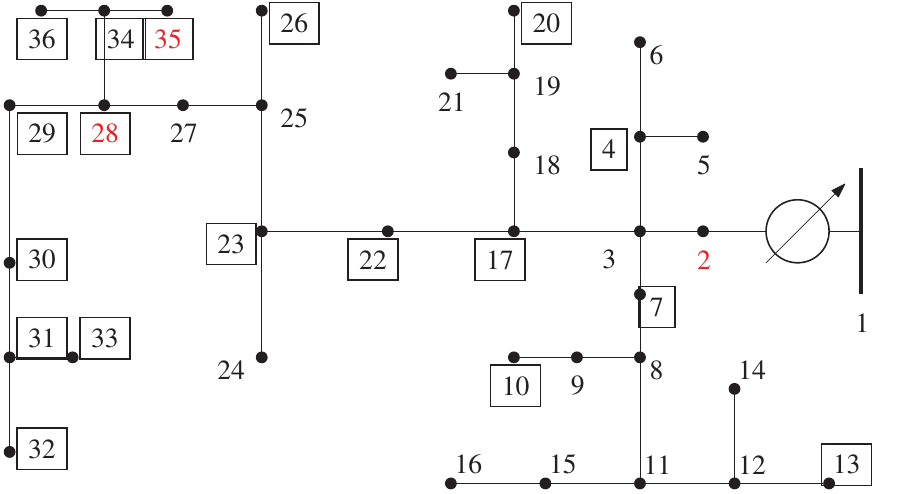}
	\caption{Modified IEEE 37-node feeder. The boxed nodes represent the location of PV systems. }
	\label{F_feeder}
\end{figure}

\begin{figure}
\centering
\includegraphics[trim = 0mm 0mm 0mm 0mm, clip, scale=0.4]{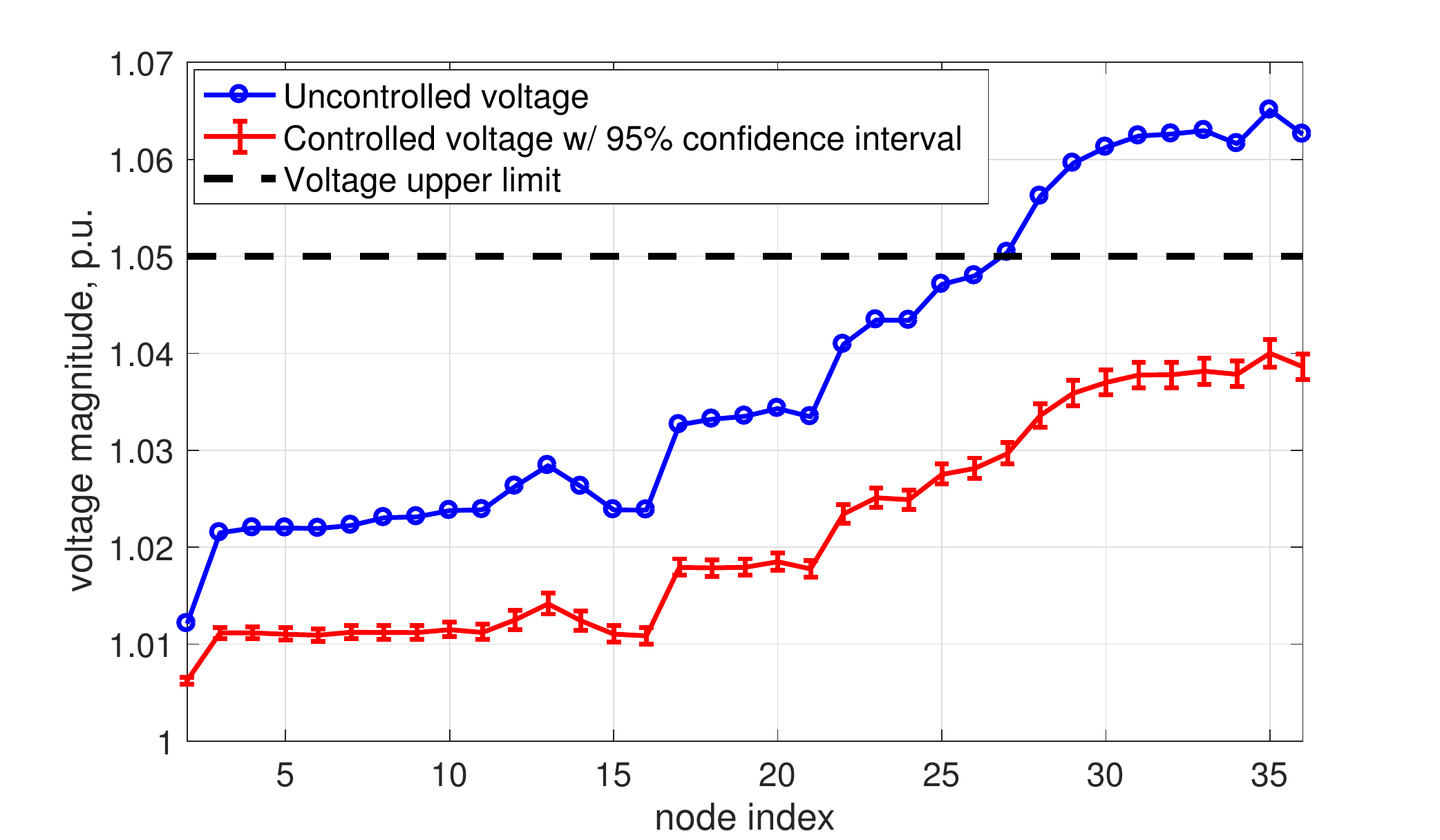}
\caption{Uncontrolled voltage and controlled voltage under scenario 2) with 95\% confidence interval.}\label{fig:controlv}
\end{figure}

Consider a modified version of the IEEE 37-node test feeder shown in Figure~\ref{F_feeder}. The modified network is obtained by considering a single-phase equivalent, and by replacing the loads specified in the original data set with real load data measured from feeders in Anatolia, California during the week of August 2012~\cite{Bank13}. 
Line impedances, shunt admittances, as well as active and reactive loads are adopted from the respective data set. It is assumed that 18 PV systems are located at nodes $4$, $7$, $10$, $13$, $17$, $20$, $22$, $23$, $26$, $28$, $29$, $30$, $31$, $32$, $33$, $34$, $35$, and $36$, constituting the set $\cD_F$, and their generation profile is simulated based on the real solar irradiance data available in~\cite{Bank13}. 
The rating of these inverters are $300$ kVA for $i = 3$, $350$ kVA for $i = 15, 16$, and $200$ kVA for the remaining inverters. With this setup, when no actions are taken to prevent overvoltages, one would obtain the voltage snapshot at noon illustrated in Figure~\ref{fig:controlv} (blue dots).

The objective functions of PV are set uniformly as $C_{i,d}(p_{i,d},q_{i,d}) =    3(p^{\textrm{av}}_{i,d} - p_{i,d})^2+q_{i,d}^2 $, where $p^{\textrm{av}}_{i,d}$ is the generated real power, in an effort to minimize the amount of real power curtailed and the amount of reactive power injected or absorbed.
We then install 15 identical TCLs at each of the following 25 nodes: 2, 5, 6, 7, 9, 10, 11, 13, 14, 16, 18, 19, 20, 21, 22, 24, 26, 27, 28, 29, 30, 32, 33, 35, and 36, totaling 375 TCLs that comprise the set $\cD_S$. We set a uniform cost function for all TCLs as $C_{i,d}(T_{i,d}^{\text{in}+})=20(T_{i,d}^{\text{in}+}-T_{i,d}^{\text{nom}})^2$, where $T_{i,d}^{\text{nom}}$ is a preferred room temperture set at  $75^{\circ}$F, and the room temperature 15 minutes later is modeled as $T_{i,d}^{\text{in}+}=T^{\text{in}}_{i,d}+0.1(T^{\text{out}}_{i,d}-T^{\text{in}}_{i,d})-0.001p_{i,d}$.
 Also, $T_{i,d}^{\text{in}+}$ should be within $[70^{\circ}\text{F},80^{\circ}\text{F}]$. For each TCL, there are two possible power rates: 0 and 4~kW.
The cost function of customer $i$ sums the cost functions of all its devices $C_i(z_i)=\sum_{d\in\cD_i}C_{i,d}(z_{i,d})$.

The voltage limits are $\overline{v}_i=1.05\!$ p.u. and $\underline{v}_i=0.95\!$ p.u., and robust voltage limits are set to $\overline{v}'_i=1.04\!$ p.u. and $\underline{v}'_i=0.96\!$ p.u. for $\forall i\!\in\!\cN$. 
Stepsize $\varepsilon=0.1$ is constant, which enables us to achieve all the results proved under diminishing stepsize. 

\begin{figure}
\centering
\includegraphics[trim = 0mm 0mm 0mm 0mm, clip, scale=0.32]{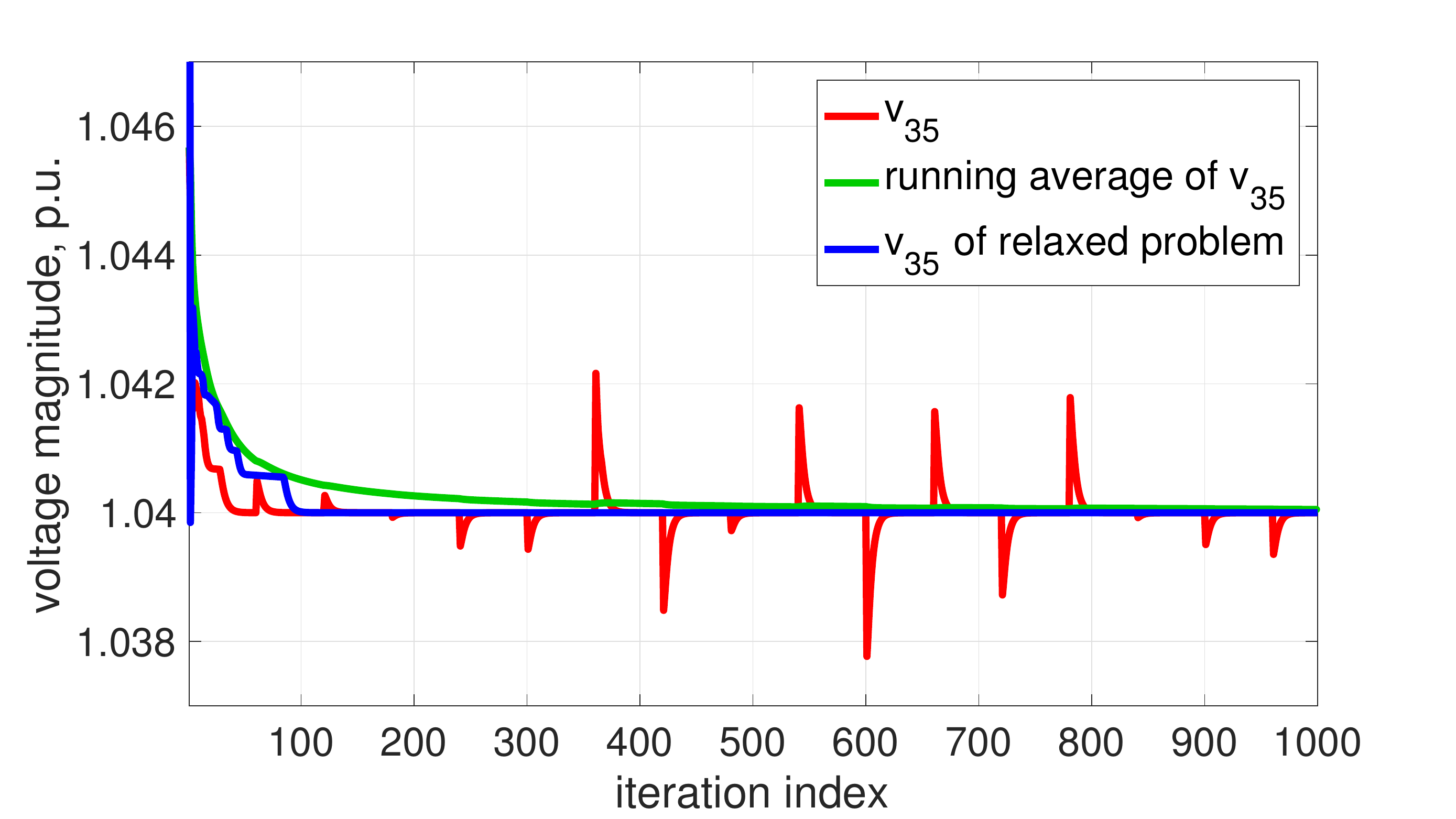}
\caption{The running average of the voltage converges to the optimal voltage of the relaxed problem.}\label{fig:runave}
\end{figure}

\subsection{Convergence and Variance}\label{sec:discretesetup}
We update PV every iteration, and TCLs every 60 iterations, with the following scenarios:
\begin{enumerate}
	\item[1)] All 15 TCLs at each node are combined to be controlled together with two power rates of 0~kW and 60~kW.
	\item[2)] 15 TCLs at each node are controlled independently.		
	\item[3)] All 15 TCLs at each node are combined to be controlled together with 16 power rates: 0~kW, 4~kW, ... , 60~kW.
\end{enumerate}

Without loss of generality, we use scenario 2) to show convergence. The results are plotted in Fig. \ref{fig:runave}. Though the resultant voltage (red line) is changing randomly, it fluctuates around the solution of the relaxed problem (blue line). Also, the running average of the fluctuating voltage (green line) approaches the solution of the relaxed problem as iteration number increases, verifying Theorem~\ref{the:convergence}--\ref{the:performance} even with constant stepsize. Moreover, in between updates of TCLs, their random consequence is absorbed by the PV's faster updates.

\begin{figure}
\centering
\includegraphics[trim = 0mm 0mm 0mm 0mm, clip, scale=0.3]{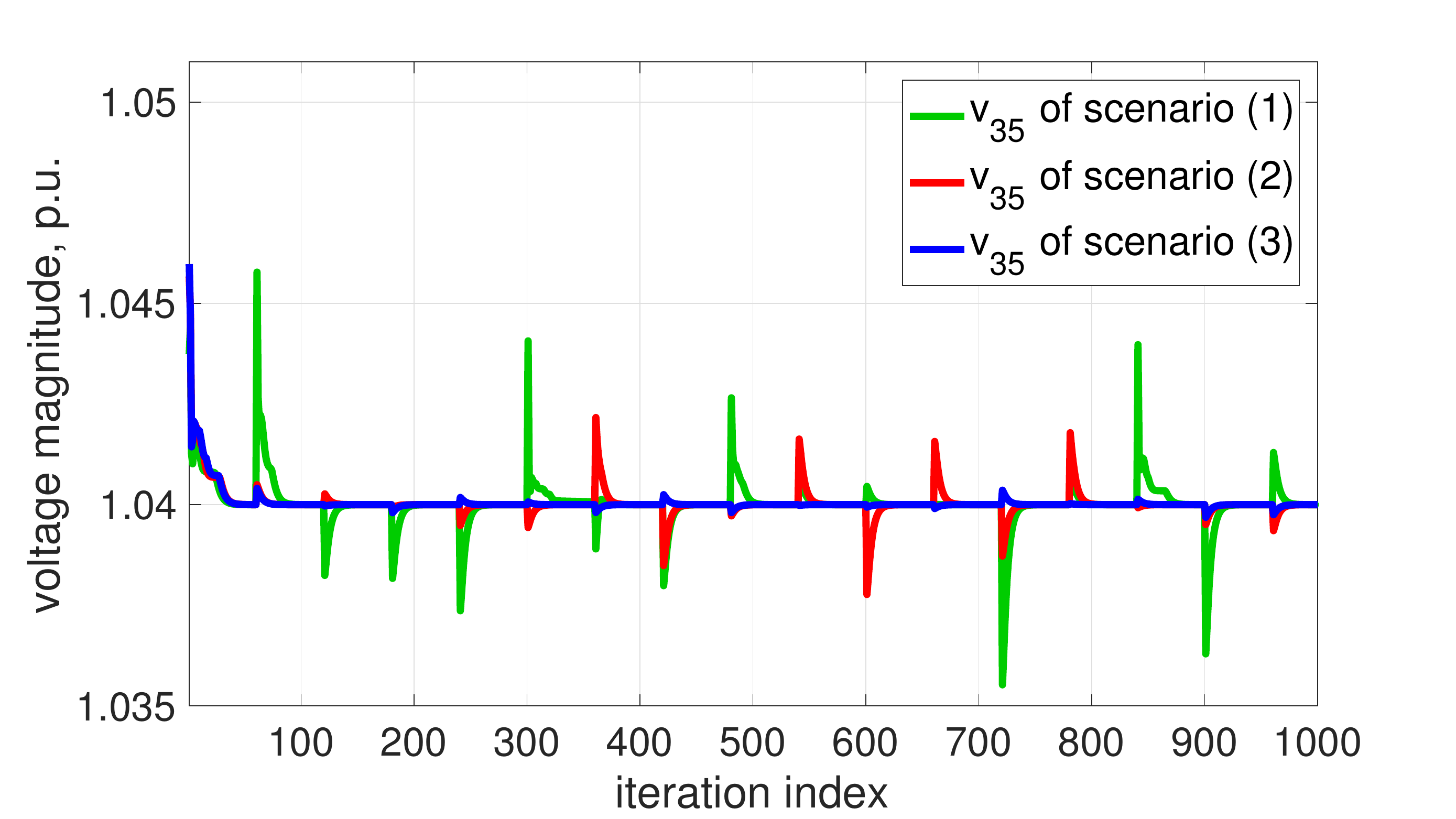}
\caption{Different voltage variance of three different scenarios.}\label{fig:ave}
\end{figure}

Next, we compare the variance of the resultant voltages among scenarios 1)--3). Based on design, these three scenarios have the same optimality on average. By Proposition~\ref{pro:var}, the voltage variance stems from the number of random variables and the control granularity. We  compare scenario 1) and 3) to illustrate that finer granularity generates less variance, and scenario 2) and 3) to show that less random variables results in less variance. The results are presented in Fig. \ref{fig:ave}, where three scenarios are marked with different colors.

\subsection{Voltage Regulation}
In this part, we use scenario 2) to compare the resultant voltages at all nodes to those without any voltage regulation. We record 25,000 random processes since convergence, and we plot the mean values of the voltages together with their 95\% confidence intervals in Fig. \ref{fig:controlv}. Because of the robust bounds as well as the small variance of the resultant voltages, the controlled voltages are all less than the original voltage upper-limit of $1.05\!$ p.u. 

\section{Conclusion}\label{sec:conclusion}
We have proposed an iterative distributed stochastic dual algorithm that allows the distribution network operator and the customers to coordinate with private information preserved to optimize the social welfare while concurrently recovering feasible power rates for discrete devices and ensuring that the voltage magnitudes are within the prescribed limits. We prove its convergence and analyze its performance. Numerical examples are provided to support the theoretical results.

\bibliographystyle{IEEEtran}
\bibliography{biblio.bib}
\newpage
\appendix
\subsection{Proof of Theorem~\ref{the:convergence}}

	The convergence of sequence $\{\mu(k)\}$ is equivalent to that of sequence $\{\|\mu(k)-\mu^*\|^2\}$, where $\mu^*$ is an optimal dual of $(\bm{\mathcal{P}_3})$. Then, at iteration $k=tM$ with some $t>0$: 
\begin{eqnarray}
	&&E\big[\|\mu(tM+1)-\mu^*\|^2| \mu(1),\ldots,\mu((t-1)M+1) \big]\nonumber\\
	&\leq&E\big[\|\mu(tM)+\varepsilon_{_{tM}} g(z(tM))-\mu^*\|^2|\mu((t-1)M+1)\big]\nonumber\\
	&\leq&E\big[\|\mu(tM)-\mu^*\|^2|\mu((t-1)M+1)\big]+\varepsilon_{_{tM}}^2G^2\nonumber\\
	&&~~+2\varepsilon_{_{tM}}(\mu(tM)-\mu^*)^{\intercal} g(z^*(tM))\nonumber\\[-2pt]
	&\leq& \|\mu((t-1)M+1)-\mu^*\|^2+\!\!\!\!\sum_{k=(t-1)M+1}^{tM}\!\!\!\!\!\!\!\varepsilon_k^2G^2\nonumber\\[-4pt]
	&&~~+\sum_{m=0}^{M-1}2\varepsilon_{_{tM-m}}(\mu(tM-m)-\mu^*)^{\intercal} g(z^*(tM-m))\nonumber\\[-4pt]
	&\leq&\|\mu((t-1)M+1)-\mu^*\|^2+\!\!\!\!\sum_{k=(t-1)M+1}^{tM}\!\!\!\!\!\!\!\varepsilon_k^2G^2\nonumber\\[-2pt]
	&&~~+2\varepsilon_{_{tM}}(h(\mu(tM))-h(\mu^*))\nonumber\\[-1pt]
	&&~~+\sum_{m=1}^{M-1}2\varepsilon_{_{tM-m}}(h_F(\mu(tM-m))-h_F(\mu^*))\nonumber\\[-4pt]
	&\leq& \|\mu((t-1)M+1)-\mu^*\|^2+\!\!\!\!\sum_{k=(t-1)M+1}^{tM}\!\!\!\!\!\!\!\varepsilon_k^2G^2,\nonumber
\end{eqnarray}
	where the first inequality comes from the non-expansiveness of projection operator, the third from repeating previous steps, the fourth from the definition of the subgradient, and the last from the definition of optimality of concave functions $h$ and $h_F$. 
	
	Because $\sum_{t=1}^{\infty}\sum_{k=(t-1)M+1}^{tM}\varepsilon_k^2G^2<\infty$, by  Lemma \ref{gladyshev} the sequence $\{\|\mu(k)-\mu^*\|^2\}$ converges to some random variable $\{\|\mu(\infty)-\mu^*\|^2\}$ almost surely, and therefore the sequence  $\{\mu(k)\}$ converges to some random vector $\mu(\infty)$ almost surely.

\subsection{Proof of Theorem~\ref{the:performance}}
	Similar to the proof of Theorem \ref{the:convergence}, we have:
	\begin{eqnarray}
	&&E[\|\mu(tM+1)-\mu^*\|^2]\nonumber\\[-6pt]
	&\leq& E\big[\|\mu((t-1)M+1)-\mu^*\|^2\big]+\!\!\!\!\sum_{k=(t-1)M+1}^{tM}\!\!\!\!\!\!\!\!\!\!\varepsilon_k^2G^2\nonumber\\[-3pt]
	&&~~+2\varepsilon_{_{tM}}(h(\mu(tM))-h(\mu^*))\nonumber\\[-1pt]
	&&~~+\sum_{m=1}^{M-1}2\varepsilon_{_{tM-m}}(h_F(\mu(tM-m))-h_F(\mu^*))\nonumber\\[-4pt]
		&\leq& E\big[\|\mu((t-1)M+1)-\mu^*\|^2\big]+\!\!\!\!\sum_{k=(t-1)M+1}^{tM}\!\!\!\!\!\!\!\!\!\!\varepsilon_k^2G^2\nonumber\\[-2pt]
	&&~~+2\varepsilon_{_{tM}}(h(\mu(tM))-h(\mu^*)).\nonumber
	\end{eqnarray}
Apply the above steps recursively to obtain:
	\begin{eqnarray}
	&&E[\|\mu(tM+1)-\mu^*\|^2] \leq E[\|\mu(1)-\mu^*\|^2]\nonumber\\
	&&\hspace{.2cm}+\sum_{k=1}^{tM} \varepsilon_k^2G^2+\sum_{\tau=1}^t 2\varepsilon_{_{\tau M}}E[h(\mu(\tau M))-h(\mu^*)].\nonumber
	\end{eqnarray}
	Because $E[\|\mu(tM+1)-\mu^*\|^2]\geq 0$, the following holds:
	\begin{eqnarray}
	&&\sum_{\tau=1}^t 2\varepsilon_{_{\tau M}}E[h(\mu^*)-h(\mu(\tau M))]\nonumber\\[-8pt]
	&\leq&E[\|\mu(1)-\mu^*\|^2]+\sum_{k=1}^{tM}\varepsilon^2_{k}G^2.\label{eq:proofstep}
	\end{eqnarray}
	By Jensen's inequality:
	\begin{eqnarray}
	E[h(\mu^*)-h(\mu(\tau M))]\geq h(\mu^*)-h(\tilde{\mu}(\tau M)).\label{eq:Jensen}
	\end{eqnarray}
	Therefore, by considering $\sum_{k=1}^{\infty}\varepsilon^2_{k}<\infty$ from (\ref{eq:stepsize}), we have from (\ref{eq:proofstep}) and (\ref{eq:Jensen}) that:
	\vspace{-.0cm}
	\begin{eqnarray}
	\lim_{t\rightarrow\infty}\sum_{\tau=1}^t 2\varepsilon_{_{\tau M}} \big(h(\mu^*)-h(\tilde{\mu}(\tau M))\big)<\infty.\label{eq:contradiction}
	\end{eqnarray}
	
	We next show $h(\mu^*)\!=\!h(\tilde{\mu}(\infty))$ by contradiction. 
	Recalling that $h(\mu^*)\!\geq\! h(\mu)$ for any feasible $\mu$, assume there exists some $e>0$ such that $h(\mu^*)\!-\!h(\tilde{\mu}(\infty))\!\geq\! e$.
	Because $\sum_{\tau=1}^{\infty}\varepsilon_{_{\tau M}}\!=\!\infty$, we must have 
	$\lim_{\tau\rightarrow\infty}\sum_{\tau=1}^t 2\varepsilon_{_{\tau M}} \big(h(\mu^*)-h(\tilde{\mu}(\tau M))\big)\!\!=\!\!\infty$, which contradicts (\ref{eq:contradiction}).
	Hence, $h(\mu^*)\!\!=\!\!h(\tilde{\mu}(\infty))$.

	Further, because $\mu$ is statistically stationary, its ensemble average and time average are identical. (\ref{eq:timeaverage}) follows.

\subsection{Proof of Proposition~\ref{pro:var}}

	The variance of $\hat{v}_i$ can be written as:
	\begin{eqnarray}
	Var(\hat{v}_i)&=&E[|\hat{v}_i-\hat{v}_i^*|^2]=E\big[\big|\sum_{j\in\cN}R_{ij}(p_j-p_j^*)\big|^2\big]\nonumber\\[-6pt]
	&\leq&\sum_{j\in\cN}R^2_{ij}\cdot E\big[\sum_{j\in\cN}(p_j-p_j^*)^2\big]\nonumber\\[-2pt]
	&\leq&\sum_{j\in\cN}R^2_{ij} \cdot \sum_{j\in\cN}\sum_{d\in\cD_{S_j}} E\big[(p_{j,d}-p_{j,d}^*)^2\big]\nonumber\\[-2pt]
	&=& \sum_{j\in\cN}R^2_{ij} \cdot  \sum_{j\in\cN}\sum_{d\in\cD_{S_j}} (p_{j,d}^*-\underline{p}_{j,d})(\overline{p}_{j,d}-p_{j,d}^*)\nonumber\\[-2pt]
	&\leq& \sum_{j\in\cN}R^2_{ij} \cdot \sum_{j\in\cN}\sum_{d\in\cD_{S_j}}{(\overline{p}_{j,d}-\underline{p}_{j,d})^2}/{4},\nonumber\\[-4pt]
	&\leq& {D_{S}}/{4}\sum_{j\in\cN}R^2_{ij}\cdot \underset{j,d}{\max} (\overline{p}_{j,d}-\underline{p}_{j,d})^2,\nonumber
	\end{eqnarray}
	where we apply Cauchy-Schwarz inequality in the first inequality, Jensen's inequality in the second, and the second equality is based on the probability distribution (\ref{eq:probability}).

\subsection{Proof of Proposition~\ref{pro:rob}}


Let $\tilde{\hat{v}}_i=E[\hat{v}_i]$. By Chebyshev's inequality \cite{stark2014probability}, given $\delta>0$, we have
\begin{eqnarray}
\mathbb{P}[|\hat{v}_i-\tilde{\hat{v}}_i|\geq \delta]\leq \frac{Var(\hat{v}_i)}{\delta^2}.\nonumber
\end{eqnarray}

Consider the upper bound first. Design the robust bound as $\overline{v}'_i \leq \overline{v}_i-\delta$, so that we must have $\tilde{\hat{v}}_i\leq \overline{v}'_i$ by Assumption~\ref{ass:sla} and Theorem~\ref{the:performance}.
And we write the probability of voltage violation as follows:
\begin{eqnarray}
&&\mathbb{P}[\hat{v}_i\geq \overline{v}_i]=\mathbb{P}[\hat{v}_i-\tilde{\hat{v}}_i\geq \overline{v}_i-\tilde{\hat{v}}_i]\nonumber\\
&\leq&\mathbb{P}[\hat{v}_i-\tilde{\hat{v}}_i\geq \overline{v}_i-\overline{v}'_i]\leq\mathbb{P}[\hat{v}_i-\tilde{\hat{v}}_i\geq \delta]\nonumber\\
&=&\frac{1}{2}\cdot\mathbb{P}[|\hat{v}_i-\tilde{\hat{v}}_i|\geq \delta]\leq \frac{Var(\hat{v}_i)}{2\delta^2}.\nonumber
\end{eqnarray}

Similar process applies to $\mathbb{P}[\hat{v}_i\leq \underline{v}_i]$.

%

\end{document}